\documentclass[letterpaper, 10 pt, conference]{ieeeconf}
\IEEEoverridecommandlockouts
\usepackage{cite}
\usepackage{amsmath,amssymb,amsfonts}

\usepackage{amsthm}
\usepackage{algorithmic}
\usepackage{import}
\usepackage{graphicx}
\usepackage{textcomp}

\usepackage{xcolor}
\usepackage{mathtools}
\let\labelindent\relax
\usepackage{enumitem}

\usepackage{tikz}
\usetikzlibrary{positioning, calc,arrows.meta}

\usepackage[absolute,overlay]{textpos} % in the preamble

\usepackage{cuted}
\usepackage{dblfloatfix} % load in the preamble

\newcommand{\Rankf}[1]{\operatorname{rank}(#1)}

\newcommand{\Dim}[1]{\mathrm{dim}(#1)}

\newcommand{\D}{\mathrm{d}}
\newcommand{\Lie}{\mathrm{L}}

\newcommand{\X}{\mathcal{X}}

\newcommand{\TX}{\mathcal{T}(\mathcal{X})}

\newcommand{\Spanbracket}[1]{\langle#1\rangle}

\newcommand{\crksub}[1]{\underset{#1}{\subset}}
\newcommand{\ddiff}[1]{d_{\operatorname{diff}}(#1)}
\newcommand{\idxsum}[1]{\vert #1 \vert}
\newcommand{\dif}[1]{d_{\operatorname{#1}}}

\newcommand{\zdot}{\dot{z}}

\newcommand{\pad}[1]{\partial_{#1}}

\newcommand{\Pc}[1]{\mathcal{P}_{#1}}
\newcommand{\Qc}[1]{\mathcal{Q}_{#1}}

\theoremstyle{plain}      % italic body
\newtheorem{theorem}{Theorem}

\newtheorem{corollary}{Corollary}

\theoremstyle{definition} % roman body
\newtheorem{definition}{Definition}
\newtheorem{example}{Example}

\theoremstyle{remark}     % roman body
\newtheorem{remark}{Remark}

\begin{document}

\title{ \textbf{A Structurally Flat Triangular Form for Three-Input Systems} 
\thanks{This work has been submitted to the IEEE for possible publication. Copyright may be transferred without notice, after which this version may no longer be accessible.}
\thanks{This research was funded in whole, or in part, by the Austrian Science Fund (FWF) P36473. For the purpose of open access, the author has applied a CC BY public copyright licence to any Author Accepted Manuscript version arising from this submission.}
}

\author{Georg Hartl, Conrad Gstöttner, and Markus Schöberl \thanks{All authors are with the Institute of Control Systems, Johannes Kepler University Linz, Altenberger Strasse 69, 4040 Linz, Austria, \mbox{E-mail:}~{\tt \{georg.hartl, conrad.gstoettner, markus.schoeberl\}@jku.at}.
}}

\maketitle

\begin{abstract}
We present a broadly applicable structurally flat triangular form 
for \mbox{$x$-flat} control-affine systems with three inputs. 
Building on recent results for the derivative structure of 
flat outputs, we define the triangular form together 
with regularity conditions that guarantee structural flatness, 
and derive necessary and sufficient conditions for a system 
with a given \mbox{$x$-flat} output to be static feedback equivalent 
to this form. Further, we present sufficient conditions under which
general \mbox{$x$-flat} three-input systems can be rendered static 
feedback equivalent to the proposed triangular form after a finite number of 
input prolongations.
\end{abstract}

% \begin{IEEEkeywords}
%     Nonlinear systems, differential flatness, feedback linearization. 
% \end{IEEEkeywords}

\section{Introduction}

Differential flatness provides a powerful 
framework for the systematic solution of feedforward and feedback 
control problems~\cite{fliess_flatness_1995, fliess_lie-backlund_1999}. A nonlinear 
control system
\begin{equation}\label{eq:f_xu_m_inputs}
    \dot x = f(x,u)
\end{equation}
with $n$ states $x$ and $m$ inputs $u$ is called flat if there 
exists an $m$-tuple of functions 
$y=\varphi(x,u,u^{(1)},\ldots,u^{(\nu)})$,
where $u^{(\nu)}$ denotes the $\nu$-th time derivative of $u$, 
such that both the state and the input can 
be expressed in terms 
of the \emph{flat output}~$y$ and finitely many of its time 
derivatives via the \emph{flat parameterization} 
$(x,u)=F(y, y^{(1)}, \ldots, y^{(r)})$.
If a flat output depends only on state variables $x$, 
we call this flat output and the corresponding system 
$x$-flat. 

However, computing a flat output for a given nonlinear multi-input 
system remains a challenging open problem, as neither practically 
verifiable necessary and sufficient conditions for flatness nor 
computationally tractable methods for computing flat outputs 
have yet been established. 

Completely solved subclasses of flat systems for which tractable 
necessary and sufficient conditions exist include \textit{static feedback 
linearizable} (SFL\footnote{We use \emph{SFL} as an abbreviation 
for \emph{static feedback linearization} and, by extension, for 
the associated notion of being \emph{statically feedback 
linearizable} and the corresponding property of \emph{static 
feedback linearizability}.}) 
systems~\cite{jakubczyk_linearization_1980}, two-input driftless 
systems~\cite{martin_feedback_1994}, general two-input systems 
that are linearizable by a two-fold prolongation of a suitably 
chosen 
control~\cite{nicolau_flatness_2016, gstottner_necessary_2023}, 
and general multi-input systems that become static feedback 
linearizable after a one-fold prolongation of a suitably chosen 
control~\cite{nicolau_flatness_2017}.

Structurally flat triangular forms have proven to be an effective 
tool for characterizing flat systems. As demonstrated, e.g., 
in~\cite{gstottner_structurally_2022, gstottner_triangular_2024, 
bououden_triangular_2011, nicolau_normal_2019, 
SilveiraFlatTriangularForm2015, schoberl_implicit_2014}, complete 
geometric characterizations of such forms yield sufficient 
conditions for flatness that require only differentiation and 
algebraic operations, and provide constructive procedures for 
computing flat outputs, which typically appear as top state 
variables in flat triangular forms.

For two-input systems, a broadly applicable \emph{general 
triangular form} ($GTF_2$) was proposed 
in~\cite{gstottner_triangular_2024}. Building on this, an 
algorithm for identifying flat-output components has been 
presented in~\cite{HartlTriangularFormsXFlat2025}. 
The $GTF_2$ encompasses the chained 
form~\cite{martin_feedback_1994}, the extended chained 
form~\cite{LiMultiinputControlaffineSystems2016, 
SilveiraFlatTriangularForm2015}, and the structurally flat 
triangular forms presented 
in~\cite{gstottner_structurally_2022, gstottner_flat_2021}. 

For systems with more than two inputs, the situation is 
considerably less developed. The canonical contact form with 
compatible drift~\cite{LiMultiinputControlaffineSystems2016} provides a 
structurally flat triangular form that has been extended 
in~\cite{hartl_flat_2025} by appending and prepending integrator 
chains. In a complementary direction, 
\cite{HartlLinearizationFlatMultiInput2026} studies when 
\mbox{$x$-flat} three-input systems can be rendered static feedback 
linearizable by a minimal number of prolongations of suitably 
chosen inputs after applying a static input transformation. 
In \cite{HartlLinearizationFlatMultiInput2026}, this property is referred to as \emph{$\varphi$-SFL via minimal prolongations}.

In this paper, we propose a general structurally flat triangular 
form for control-affine three-input systems of the form
\begin{equation}\label{eq:f_xu_affine_three_inputs}
    \dot x = f(x) + g_1(x)u^1 + g_2(x)u^2 + g_3(x)u^3,
\end{equation}
that admit an \mbox{$x$-flat} output
\begin{equation}\label{eq:x_flat_output_three_inputs}
    y = (\varphi^1(x), \varphi^2(x), \varphi^3(x)) \, .
\end{equation}
The proposed form encompasses several previously characterized 
structurally flat triangular forms. Our contribution is as follows:
\begin{enumerate}[label=\Alph*)]
    \item We introduce the general structurally flat triangular 
    form for \mbox{$x$-flat} control-affine three-input systems and 
    present necessary and sufficient conditions for a system 
    with a given \mbox{$x$-flat} output to be static feedback equivalent 
    (SFE) to this form. 
    \item We prove that every \mbox{$x$-flat} 
    system~\eqref{eq:f_xu_affine_three_inputs} that is $\varphi$-SFL
    via minimal prolongations can be transformed into 
    the general triangular form after a finite 
    number of input prolongations.
\end{enumerate}

Our paper is structured as follows: 
Section~\ref{sec:notation} introduces notation and terminology. 
Section~\ref{sec:Known_Results} recalls the derivative structure 
of three-input \mbox{$x$-flat} outputs. Our main results are presented in 
Section~\ref{sec:main_results}. Illustrating examples are given in 
Section~\ref{sec:examples}. Finally, 
Section~\ref{sec:conclusion} gives a conclusion and an outlook 
for further research. Remaining proofs are presented in the 
appendix.

%% ================================================================
\section{Notation and Terminology}
\label{sec:notation}

Throughout, the Einstein summation convention and the tensor notation
are used omitting the index range when clear from the context.
For a smooth manifold 
$\mathcal{X}$ with local 
coordinates~\mbox{$x=(x^1, \ldots, x^n)$},
the \textit{tangent} and \textit{cotangent bundle} on $\X$ are denoted by
$\TX$ and $\mathcal{T}^*(\X)$. Given an $m$-tuple of 
functions~\mbox{$h=(h^1, \ldots, h^m):\mathcal{X} 
\rightarrow \mathbb{R}^m$}, we write its $m \times n$ 
Jacobian matrix as $\partial_x h$. In particular, 
$\partial_{x^i}h^j$ represents the partial derivative of $h^j$ 
with respect to $x^i$. The notation $\D h$ represents the 
differentials~\mbox{$(\D h^1, \ldots, \D h^m)$}. Given a 
set of one-forms $\omega = (\omega^1, \ldots, \omega^m)$ 
that defines a \emph{codistribution} 
$\Pc{} = \Spanbracket{\omega^1, \ldots, \omega^m} \subset \mathcal{T}^*(\X)$, 
$\Spanbracket{\omega}$ is used to 
indicate the set of all possible linear combinations 
of~\mbox{$(\omega^1, \ldots, \omega^m)$}. The notation \(\Pc{1}\underset{k}{\subset}\Pc{2}\) indicates that the corank of \(\Pc{1}\) in \(\Pc{2}\) is \(k\). Given a \textit{vector field} $v \in \TX$, we denote the 
Lie derivative of a scalar function $h^j$ along $v$ by $\Lie_vh^j$. 
 
To represent time derivatives we use subscripts in square 
brackets. For example, $y^j_{[\alpha]}$ represents the 
$\alpha$-th time derivative of the $j$-th element within the 
$m$-tuple $y$, 
and~\mbox{$y_{[\alpha]} = (y^1_{[\alpha]}, \dots, 
y^m_{[\alpha]})$} indicates the $\alpha$-th time derivative for 
each component of $y$. Capitalized multi-indices are used to represent time 
derivatives of different orders of each element in a tuple. Given 
two multi-indices \mbox{$A = (a^1, \dots, a^m)$} 
and~\mbox{$B=(b^1, \dots, b^m)$} with $a^j \leq b^j$, denoted 
by $A \leq B$, we can find concise representations of time 
derivatives of different orders, such as 
\mbox{$y_{[A]} = (y^1_{[a^1]}, \dots,y^m_{[a^m]})$}, 
\mbox{$y_{[0, A]} = ( y^1_{[0,a^1]}, \dots,y^m_{[0,a^m]} )$}, 
and 
\mbox{$y_{[A, B]} = ( y^1_{[a^1,b^1]}, \dots,y^m_{[a^m,b^m]} )$}. 
Then \( y^j_{[a^j, b^j]} \) denotes successive time derivatives 
of $y^j$ given by~\mbox{\( y^j_{[a^j, b^j]} = ( y^j_{[a^j]}, \dots, 
y^j_{[b^j]}) \)}. If~\mbox{$a^j > b^j$}, then 
\mbox{$y^j_{[a^j, b^j]}$} is empty. We express the 
component-wise addition and subtraction of multi-indices 
as~\mbox{$A \pm B = ( a^1\pm b^1, \dots, a^m \pm b^m )$}. 
The summation over the indices is given by 
$\idxsum{A} = \sum_{j=1}^m a^j$. 

We use the overline to represent successive state variables 
of a subsystem. That is, 
$\overline{x}^i_j = (x^1_j, \ldots, x^i_j)$, where the 
superscript $i$ indicates the range and the subscript $j$ 
distinguishes different subsystems. Given a subsystem 
consisting of $k$ states, for the complete collection 
we drop the overline and the superscript, and simply 
write $x_j = (x^1_j, \ldots, x^k_j) (= \overline{x}^k_j)$.

Throughout, when coefficient functions appear in the same 
expression, e.g., 
$a^{k^2}_2(z_1, \overline{z}_2^{k^2+1}, 
\overline{z}_3^{k^3+\delta+1}) 
+ b^{k^2}_{2,1}(\ldots)\hat u^1$, 
the abbreviation $(\ldots)$ indicates that 
$b^{k^2}_{2,1}$ may depend on the same variables as 
$a^{k^2}_2$.

A nonlinear system~\eqref{eq:f_xu_m_inputs} is 
\emph{static feedback equivalent} (SFE) to another system 
$\zdot = g(z, \hat u)$, if there exist a state transformation 
$z=\Phi_z(x)$ and a static input transformation 
$\hat u=\Phi_{\hat u}(x,u)$ such that the system dynamics satisfy
\vspace{-0.5ex}
\begin{equation*}
    \left(f^i(x,u)\partial_{x^i}\Phi^j_z(x)\right)
    \circ(\Phi_x(z), \Phi_{u}(z,\hat u))=g^j(z, \hat u),
    \vspace{-0.5ex}
\end{equation*}
where $(x, u)=(\Phi_x(z), \Phi_{u}(z,\hat u))$ denotes the inverse of 
$(z,\hat u)=(\Phi_z(x),\Phi_{\hat u}(x,u))$. Unless stated otherwise, any given transformation is considered 
invertible. We assume that all functions, vector fields, and covector 
fields are smooth, and that all codistributions have locally 
constant rank.

%%%%%%%%%%%%%%%%%%%%%%%%%%%%%%%%%%%%%%%%%%%%%%%%%%%%%%%%%%%%%%%%%%
%%%%%%%%%%%%%%%%%%%%%%%%%%%%%%%%%%%%%%%%%%%%%%%%%%%%%%%%%%%%%%%%%%
\begin{figure*}[!b]
\hrulefill
\normalsize
\begin{equation}\label{eq:deriv_form_case_rank1}
    \tag{$DS_3$}
    \begin{alignedat}{3}
        \varphi^{1}_{[0,k^1-1]} & = \varphi^{1}_{[0,k^1-1]}(x) 
        & \hspace{1em} & \varphi^{2}_{[0,k^2-1]}(x) 
        & \hspace{1em} & \varphi^{3}_{[0,k^3-1]}(x) \\
        \varphi^{1}_{[k^1]} & = \hat u^1 
        &  &  \varphi^{2}_{[k^2]}( x, \hat u^1 ) 
        &  &  \varphi^{3}_{[k^3]}( x, \hat u^1 ) \\[-0.5ex]
        & \hspace{0.5em} \vdots 
        && \hspace{0.5em} \vdots 
        && \hspace{0.5em} \vdots \\[-0.5ex]
        & &  & \varphi^2_{[p^2]}(x, \hat u^1_{[0,p^2-k^2]}, u^2) 
        &  & \varphi^{3}_{[p^3]}( x, \hat u^1_{[0,p^3-k^3]}, u^2) \\[-0.5ex]
        & \hspace{0.5em} \vdots 
        && \hspace{0.5em} \vdots 
        && \hspace{0.5em} \vdots \\[-0.5ex]
        & & & \varphi^2_{[p^2+s]}(x, \hat u^1_{[0,p^2-k^2+s]}, u^2_{[0,s]}, u^3) 
        & & \varphi^{3}_{[p^3+s]}( x, \hat u^1_{[0,p^3-k^3+s]}, u^2_{[0,s]}, u^3) \\[-0.5ex]
        & \hspace{0.5em} \vdots 
        && \hspace{0.5em} \vdots 
        && \hspace{0.5em} \vdots \\[-0.5ex]
        \varphi^{1}_{[r^1-1]} & = \hat u^1_{[\dif{max}-1]} 
        & & \varphi^2_{[r^2-1]}(x, \hat u^1_{[0,\dif{min}-1]}, u^2_{[0,d_{rp}-1]}, u^3_{[0,d_{rp}-s-1]}) 
        & & \varphi^{3}_{[r^3-1]}( x, \hat u^1_{[0,\dif{max}-1]}, u^2_{[0,d_{rp}-1]}, u^3_{[0,d_{rp}-s-1]} )\\
        \varphi^{1}_{[r^1]} & = \hat u^1_{[\dif{max}]} 
        & & \varphi^2_{[r^2]}(x, \hat u^1_{[0,\dif{min}]}, u^2_{[0,d_{rp}]}, u^3_{[0,d_{rp}-s]})
        & & \varphi^{3}_{[r^3]}( x, \hat u^1_{[0,\dif{max}]}, u^2_{[0,d_{rp}]}, u^3_{[0,d_{rp}-s]} )\\
    \end{alignedat}
\end{equation}
% \hrulefill
\end{figure*}
%%%%%%%%%%%%%%%%%%%%%%%%%%%%%%%%%%%%%%%%%%%%%%%%%%%%%%%%%%%%%%%%%%
%%%%%%%%%%%%%%%%%%%%%%%%%%%%%%%%%%%%%%%%%%%%%%%%%%%%%%%%%%%%%%%%%%

\subsection{Differential Flatness}
 
Within a differential geometric setting as, e.g., 
in~\cite{kolar_properties_2016}, we work on an extended 
state-input manifold $\X \times \mathcal{U}_{[0,l_u]}$ with 
local coordinates $(x,u_{[0,l_u]})$, where the integer $l_u$ is 
chosen large enough to ensure that all time derivatives of 
relevant functions along trajectories of \eqref{eq:f_xu_m_inputs}
can be computed via Lie derivatives along the 
vector field \vspace{-1ex}
\begin{equation}\label{eq:extended_vector_field}
    f_u=f^i(x,u)\partial_{x^i}
    +\sum_{\alpha=0}^{l_u-1}u^j_{[\alpha+1]}
    \partial_{u^j_{[\alpha]}}.
\end{equation}
 
\begin{definition}\label{def:diff_flatness}
A system of the form~\eqref{eq:f_xu_m_inputs} is differentially 
flat if there exist $m$ smooth functions
\begin{equation}\label{eq:flat_output}
    y=\varphi(x,u_{[0,Q]}),
\end{equation}
on $\X \times \mathcal{U}_{[0,l_u]}$ that permit a local 
parametrization
\begin{equation}
    \label{eq:flat_parametrization}
    \begin{alignedat}{2}
        x^i & = F_x^i(y_{[0,R-1]}), & \hspace{3em} 
        i&=1,\ldots,n, \\
        u^j & = F_u^j(y_{[0,R]}), & j&=1,\ldots,m,
    \end{alignedat}
\end{equation}
where $F^i_x$ and $F^j_u$ are smooth functions, and the 
multi-index $R=(r^1,\ldots,r^m)$ specifies the highest 
derivative-orders of $y$ appearing 
in~\eqref{eq:flat_parametrization}. The $m$-tuple of 
functions~\eqref{eq:flat_output}, where $Q$ defines the 
highest-derivative orders of $u$ therein, is called a flat output 
of the system.
\end{definition}
 
As shown in~\cite{kolar_properties_2016}, the differentials 
$\D \varphi, \D \varphi_{[1]}, \ldots, \D \varphi_{[\beta]}$ 
are linearly independent for arbitrary differentiation 
order $\beta$. This guarantees local uniqueness of the 
parameterization~\eqref{eq:flat_parametrization} and of the 
multi-index $R$. The 
mapping \mbox{$(F_x,F_u):\mathbb{R}^{\idxsum{R}+m} \to 
\mathbb{R}^{n+m}$} is a submersion, whose corank defines the 
\emph{differential difference} 
\mbox{$\dif{diff}(\varphi) = \idxsum{R}-n$}. Flatness further 
implies 
\begin{equation}\label{eq:imp_exterior_deriv}
    \Spanbracket{\D x}  \subset 
    \Spanbracket{\D \varphi_{[0,R-1]}}, 
    \qquad 
    \Spanbracket{\D u}  \subset 
    \Spanbracket{\D \varphi_{[0,R]}} . 
\end{equation}

% A system~\eqref{eq:f_xu_m_inputs} that admits a flat output of 
% the form \mbox{$y=\varphi(x,u)$} is called $(x,u)$-flat. 
Each component $\varphi^j(x)$ of an $x$-flat output is characterized 
by its relative degree $k^j$, defined via the vector 
field~\eqref{eq:extended_vector_field} as
\begin{equation*}
    \Lie_{f_u}^{k^j-1}\varphi^j=\varphi^j_{[k^j-1]}(x),
    \quad 
    \Lie_{f_u}^{k^j}\varphi^j=\varphi^j_{[k^j]}(x,u) \, .
\end{equation*}
% If a component $\varphi^j(x,u)$ already depends explicitly on 
% an input, the corresponding relative degree is $k^j=0$.

The following definition formalizes the linearizability of flat systems
via static feedback after prolongations of suitably 
chosen inputs.

\begin{definition}\label{def:SFL}
    Consider a system~\eqref{eq:f_xu_m_inputs} with a flat 
    output $\varphi$ of the form~\eqref{eq:flat_output}. The 
    given system is \textit{$\varphi$-SFL via prolongations} if 
    there exist an invertible static input transformation 
    $(\hat u_1, \hat u_2) = \Phi_{\hat u}(x,u)$, where 
    $\Dim{\hat u_1}=m_1$, and a 
    multi-index~\mbox{$D=(d^1,\ldots,d^{m_1})$}, with each 
    $d^j >0$, such that 
    \begin{equation*}
        \begin{aligned}
            \dot x = \hat{f}(x,\hat u_1, \hat u_2),\; 
            \dot {\hat u}_1 = \hat u_{1,[1]},\; \ldots, \; 
            \dot {\hat u}_{1,[D-1]} = \hat u_{1,[D]}
        \end{aligned}
    \end{equation*}
    is SFL with $\varphi$ as linearizing output. If 
    $\idxsum{D} = d_{\mathrm{diff}}(\varphi)$ holds 
    additionally, the system is called 
    $\varphi$\textit{-SFL via minimal prolongations}. 
\end{definition}
A detailed analysis of $\varphi$-SFL via minimal prolongations 
for $x$-flat three-input systems is given 
in~\cite{HartlLinearizationFlatMultiInput2026}.
 
%% ================================================================
\section{Known Results}
\label{sec:Known_Results}
 
In this section, we recall important structural properties of \mbox{$x$-flat} 
three-input systems from \cite{HartlLinearizationFlatMultiInput2026} 
that form the basis for the definition and characterization of the proposed 
triangular form. The following theorem, based 
on~\cite[Theorem~2]{HartlLinearizationFlatMultiInput2026}, 
describes the structure of the time derivatives of an \mbox{$x$-flat} 
output after applying a suitable input transformation. 

\begin{theorem}\label{thm:rank1_case}
    Consider a system of the 
    form~\eqref{eq:f_xu_affine_three_inputs} that admits an 
    \mbox{$x$-flat} output~\eqref{eq:x_flat_output_three_inputs}, 
    characterized by the relative degrees $K$ and the multi-index 
    $R$. First, the flat-output components can always be 
    rearranged such that
    \begin{equation}\label{eq:rearrangement_case_rank1}
        \dif{max} \coloneq r^1-k^1 = r^3-k^3 \geq r^2-k^2 \eqqcolon \dif{min}.
    \end{equation}
    
    Then, after a possible relabeling of the input components, 
    let $\hat u^1 = \varphi^1_{[k^1]}(x,u)$ replace $u^1$. Under 
    this input transformation, the following properties hold: 
    
    a) The derivatives $\varphi_{[0,R]}$ take the 
    form~\eqref{eq:deriv_form_case_rank1}, where $p^2$ and $p^3$ 
    denote the smallest derivative orders such that 
    $\varphi^j_{[p^j]}$, \mbox{$j\in\{2,3\}$}, explicitly depend 
    on $u^2$, and $s \geq 0$ is the smallest integer such that 
    $\varphi^j_{[p^j+s]}$, $j\in\{2,3\}$, explicitly depend on 
    $u^3$. Furthermore,
    \begin{equation*}
        r^2-p^2 = r^3-p^3 \eqqcolon d_{rp}, \hspace{1em} 
        \text{and} \hspace{1em} p^3-k^3 \geq p^2-k^2 \, .
    \end{equation*}
    
    b) The state dimension and $\ddiff{\varphi}$ are given by
    $\ddiff{\varphi} = \dif{max}+d_{rp}$, and 
    $n = k^1+p^2+r^3$.

    c) The system is \mbox{$\varphi$-SFL} via minimal prolongations 
    with an input transformation of the form
    \begin{equation*}\label{eq:input_trf_case_rank1}
        \hat u_1 = (\hat u^1,\hat u^2) =
        \bigl(\varphi^1_{[k^1]}(x,u),
        \phi_{\hat u^2}(x,u)\bigr),
    \end{equation*}
    replacing $(u^1,u^2)$, and $D=(\dif{max}, d_{rp})$, 
    if and only if there exists a transformation 
    \mbox{$\hat u^2 = \phi_{\hat u^2}(x, u)$} 
    such that \mbox{$s = d_{rp}$}.
\end{theorem}

\begin{remark}
    For $p^2=k^2$, i.e., for \mbox{$\Rankf{\pad{u}(\varphi_{[K]})}=2$}, 
    a possible input transformation satisfying item c) of Theorem~\ref{thm:rank1_case}
    is given by \mbox{$\phi_{\hat u^2} = \varphi^2_{[k^2]}$}.
\end{remark}

%%%%%%%%%%%%%%%%%%%%%%%%%%% GTF 3 %%%%%%%%%%%%%%%%%%%%%%%%%%%%%%%%%%%%%
\begin{figure*}[!b]
\hrulefill
\normalsize
\vspace{-1.3ex}
\begin{equation}
    \tag{$GTF_3$}
    \label{eq:GTF_3_inputs}
    \begin{aligned}[t]
        \dot z^1_1 & = z^2_1 \\[-1.5ex]
        & \hspace{0.5em} \vdots \\[-1.5ex]
        \dot z^{k^1-1}_1 & = z^{k^1}_1 \\
        \dot z^{k^1}_1 & = \hat u^1 \\
    \end{aligned}
    \quad
    \begin{aligned}[t]
        \dot z^1_2 & = z^2_2 \\[-1.5ex]
        & \hspace{0.5em} \vdots \\[-1.5ex]
        \dot z^{k^2-1}_2 & = z^{k^2}_2 \\[-1ex]
        \\
        \\
        \dot z^{k^2}_2 & = 
        a^{k^2}_2( z_1, \overline{z}_2^{k^2+1}, 
        \overline{z}_3^{k^3 + \delta + 1} )
        + b^{k^2}_{2,1}( \ldots )\hat u^1 \\[-1.5ex]
        & \hspace{0.5em} \vdots \\[-1.5ex]
        \dot z^{p^2-1}_2 & = 
        a^{p^2-1}_2( z_1,z_2, \overline{z}_3^{p^3} )
        + b^{p^2-1}_{2,1}( \ldots )\hat u^1 \\
        % & \hspace{2em}
        % + b^{p^2-1}_{2,1}( z_1, \overline{z}_2^{p^2}, \overline{z}_3^{p^3} )\hat u^1 \\ \\[-1ex]
        \dot z^{p^2}_2 & = \hat u^2
    \end{aligned}
    \quad
    \begin{aligned}[t]
        \dot z^1_3 & = z^2_3 \\[-1.5ex]
        & \hspace{0.5em} \vdots \\[-1.5ex]
        \dot z^{k^3-1}_3 & = z^{k^3}_3 \\
        \dot z^{k^3}_3 & = a^{k^3}_3( z_1, \overline{z}_2^{k^2-\delta+1},
        \overline{z}_3^{k^3+1} ) 
        + b^{k^3}_{3,1}( \ldots )\hat u^1 \\[-1.5ex]
        & \hspace{0.5em} \vdots \\[-1.5ex]
        \dot z^{k^3+\delta}_3 & = 
        a^{k^3+\delta}_3( z_1, 
        \overline{z}_2^{k^2+1}, 
        \overline{z}_3^{k^3 + \delta + 1} ) +
        b^{k^3+\delta}_{3,1}( \ldots )\hat u^1 \\[-1.5ex]
        & \hspace{0.5em} \vdots \\[-1.5ex]
        \dot z^{p^3\!-\!1}_3 & = 
        a^{p^3\!-\!1}_3( z_1,  z_2, \overline{z}_3^{p^3}  ) 
        + b^{p^3\!-\!1}_{3,1}( \ldots )\hat u^1 \\
        \dot z^{p^3}_3 & = 
        a^{p^3}_3( z_1,  z_2, \overline{z}_3^{p^3+1} ) 
        + b^{p^3}_{3,1}( \ldots )\hat u^1 + b^{p^3}_{3,2}( \ldots )\hat u^2 \\[-1.5ex]
        & \hspace{0.5em} \vdots \\[-1.5ex]
        \dot z^{r^3-1}_3 & = a^{r^3-1}_3(z) 
        + b^{r^3-1}_{3,1}(z)\hat u^1 
        + b^{r^3-1}_{3,2}(z)\hat u^2 \\
        \dot z^{r^3}_3 & = \hat u^3
    \end{aligned}\vspace{-1ex}
\end{equation}
% \hrulefill
\end{figure*}
\section{Main Results}
\label{sec:main_results}

This section introduces the general triangular form for 
\mbox{$x$-flat} three-input systems and provides a geometric characterization 
based on a given flat output $\varphi(x)$.
Subsequently, we establish a relation 
between $\varphi$-SFL via minimal prolongations and the proposed triangular form.

\subsection{ A General Flat Triangular Form for Three-Input Systems }

Consider a control-affine system~\eqref{eq:f_xu_affine_three_inputs} with
a flat output~\eqref{eq:x_flat_output_three_inputs} and let its components 
be arranged such that~\eqref{eq:rearrangement_case_rank1} holds.
According to Theorem~\ref{thm:rank1_case}, applying the input 
transformation $\hat u^1 = \varphi^1_{[k^1]}(x,u)$ yields time derivatives 
$\varphi_{[0,R]}$ of the form \eqref{eq:deriv_form_case_rank1}.
For notational convenience,
we introduce 
\[
\delta = \dif{max} - \dif{min} \, .
\]

From~\eqref{eq:rearrangement_case_rank1} and item~a) of Theorem~\ref{thm:rank1_case} 
it follows that
\begin{equation*}
    p^3 - k^3 \geq p^2 - k^2, \hspace{1em} \text{ and } \hspace{1em} 
    (p^3 - k^3) - (p^2 - k^2) = \delta.
\end{equation*}
We further define the multi-index $P = (p^1, p^2, p^3)$, where 
$p^1$ is chosen such that $p^1 - k^1 = p^3 - k^3$. 

Based on the derivative structure established in 
Theorem~\ref{thm:rank1_case}, we establish the general 
structurally flat triangular form for three-input systems.

\begin{definition}\label{def:GTF3}
Consider a three-input system~\eqref{eq:GTF_3_inputs} with 
$n=k^1+p^2+r^3$ states 
\mbox{$z=(\overline{z}^{k^1}_1,\overline{z}^{p^2}_2,
\overline{z}^{r^3}_3)$}. 
Suppose that 
$b^{k^2}_{2,1}\neq 0$, $b^{k^3}_{3,1} \neq 0$, and that the 
regularity conditions
\begin{equation}\label{eq:reg_cond_1}
    \pad{z^{k^3+i+1}_3}(a^{k^3+i}_3
    +b^{k^3+i}_{3,1} \hat u^1) \neq 0
\end{equation}
for $0 \leq i \leq \delta-1$, 
\begin{equation}\label{eq:reg_cond_2}
    \operatorname{det}\left(
    \pad{(z^{k^2+i+1}_{2}, z^{k^3+\delta+i+1}_{3})}
    \begin{pmatrix}
        a^{k^2+i}_2+b^{ k^2+i }_{ 2,1 } \hat u^1 \\
        a^{k^3+\delta+i}_3+b^{ k^3+\delta+i }_{ 3,1 } \hat u^1 
    \end{pmatrix}
    \right) \neq 0
\end{equation}   
for $0 \leq i \leq p^2-k^2-1$, and
\begin{equation}\label{eq:reg_cond_3}
    \pad{z^{k^3+i+1}_3}(a^{k^3+i}_3
    +b^{k^3+i}_{3,1} \hat u^1
    +b^{k^3+i}_{3,2} \hat u^2) \neq 0
\end{equation}
for $(p^3-k^3) \leq i \leq r^3 - k^3-1$ hold. Then we refer 
to~\eqref{eq:GTF_3_inputs} as the \emph{general structurally 
flat triangular form} for \mbox{$x$-flat} three-input systems with the 
corresponding flat output
\begin{equation}\label{eq:flat_output_GTF3}
    \varphi = (z^1_1, z^1_2, z^1_3).
\end{equation} 
\end{definition}
\smallskip
The triangular form~\eqref{eq:GTF_3_inputs} consists of three 
coupled subsystems: the $z_1$-subsystem given by an integrator 
chain of dimension $k^1$ with the input $\hat u^1$, 
the $z_2$-subsystem of dimension $p^2$ where 
$\hat u^1$ and $\hat u^2$ act as inputs, and the \mbox{$z_3$-subsystem} of 
dimension $r^3$ in which all three inputs enter in a triangular 
way. 

Given a system of the form~\eqref{eq:GTF_3_inputs}, 
the regularity conditions ensure that 
each level of the dynamics uniquely determines the parameterization 
of the succeeding set of state variables in terms of the flat output
\eqref{eq:flat_output_GTF3} and finitely many of its time derivatives.
The flat parameterization~\eqref{eq:flat_parametrization} 
can then be computed from top to bottom as follows:
\begin{enumerate}[label=\arabic*)]
    \item The integrator chains yield 
    $\overline{z}_j^{k^j} = y^j_{[0,k^j-1]}$ for 
    \mbox{$j \in \{1,2,3\}$}, and $\hat u^1 = y^1_{[k^1]}$.
    \item The states $z_3^{k^3+1}, \ldots, z_3^{k^3+\delta}$ 
    are determined from 
    $\dot z_3^{k^3}, \ldots, \dot z_3^{k^3+\delta-1}$, 
    guaranteed by~\eqref{eq:reg_cond_1}.
    \item By~\eqref{eq:reg_cond_2}, 
    the states $z_2^{k^2+1}, \ldots, z_2^{p^2}$ and 
    $z_3^{k^3+\delta+1}, \ldots, z_3^{p^3}$ are determined 
    from $\dot z_2^{k^2}, \ldots, \dot z_2^{p^2-1}$ and 
    $\dot z_3^{k^3+\delta}, \ldots, \dot z_3^{p^3-1}$. The input $\hat u^2$ 
    follows from \mbox{$\dot z_2^{p^2} = \hat u^2$}.
    \item The remaining states 
    $z_3^{p^3+1}, \ldots, z_3^{r^3}$ are determined from 
    $\dot z_3^{p^3}, \ldots, \dot z_3^{r^3-1}$, guaranteed 
    by~\eqref{eq:reg_cond_3}. The input $\hat u^3$ follows 
    from $\dot z_3^{r^3} = \hat u^3$.
\end{enumerate}

In the following, we establish necessary and sufficient
geometric conditions for static feedback equivalence 
to~\eqref{eq:GTF_3_inputs} formulated in terms of a sequence
of integrable codistributions associated with a given $x$-flat
output. To define a sequence of 
codistributions whose integrability characterizes 
the triangular form~\eqref{eq:GTF_3_inputs}, we first introduce
the multi-index 
\begin{equation}\label{eq:multi_index_A}
    A(i) = (k^1-1+i, \, k^2-1+\max(i-\delta,0), \, k^3-1+i) \, 
\end{equation}
for $i=0, \ldots, \dif{max}$. Note that \mbox{$A(0)=K-1$}, 
\mbox{$A(p^3-k^3)=P-1$}, and \mbox{$A(\dif{max})=R-1$}. 
% Using the multi-index \eqref{eq:multi_index_A}, 
For any system with a given $x$-flat output, there exists 
a corresponding sequence of codistributions
\begin{equation*}
    \Pc{A(i)} = \Spanbracket{\D \varphi_{[0,A(i)]}} \, , 
    \quad i=0,\ldots \dif{max}\, , 
\end{equation*}
taking the form
\begin{subequations}\label{eq:P_sequ}
\begin{equation}\label{eq:P_sequ_to_delta}
\begin{alignedat}{1}
    \Pc{A(0)} &= \Spanbracket{ \D \varphi_{[0,K-1]} }, \\
    % \Pc{A(1)} &= \Spanbracket{ \D \varphi^1_{[0,k^1]}, \D \varphi^2_{[0,k^2-1]}, \D \varphi^3_{[0,k^3]} } \\[-1.5ex]
    & \hspace{0.5em} \vdots \\[-0.6ex]
    \Pc{A(\delta)} &= \Spanbracket{ \D \varphi^1_{[0,(k^1 - 1) + \delta]}, \D \varphi^2_{[0, k^2-1]}, \D \varphi^3_{[0,(k^3 - 1) + \delta]} }, \\
\end{alignedat}
\end{equation}
where $\D \varphi^2_{[k^2]}$ does not appear, 
since $A(i)$ only increments the 
first and third components for $i \leq \delta$, and
\begin{equation}\label{eq:P_sequ_delta_plus}
\begin{alignedat}{1}
    \Pc{A(\delta+1)} &= \Spanbracket{ \D \varphi^1_{[0,k^1 + \delta]}, \D \varphi^2_{[0,k^2]}, \D \varphi^3_{[0,k^3 + \delta]} }, \\[-1.5ex]
    & \hspace{0.5em} \vdots \\[-0.6ex]
    \Pc{A(p^3-k^3+1)} &= \Spanbracket{ \D \varphi^1_{[0,p^1]}, \D \varphi^2_{[0,p^2]}, \D \varphi^3_{[0,p^3]} }, \\[-1.5ex]
    % \Pc{A(p^3-k^3+1)} &= \Spand{ \D \varphi^1_{[0,p^1]}, \D \varphi^2_{[0,p^2]}, \D \varphi^3_{[0,p^3]} } \\
    % \Pc{A(p^3-k^3+2)} &= \Spand{ \D \varphi^1_{[0,p^1+1]}, \D \varphi^2_{[0,p^2+1]}, \D \varphi^3_{[0,p^3+1]} } \\
    & \hspace{0.5em} \vdots \\[-0.6ex]
    \Pc{A(\dif{max})} &= \Spanbracket{ \D \varphi^1_{[0,r^1-1]}, \D \varphi^2_{[0,r^2-1]}, \D \varphi^3_{[0,r^3-1]} } \, .
\end{alignedat}
\end{equation}
\end{subequations}

Intersecting each codistribution of~\eqref{eq:P_sequ} 
with the span of the state differentials, i.e.,
\begin{equation*}
    \Qc{A(i)} = \Pc{A(i)} \cap \Spanbracket{\D x}, \qquad 
    i = 0, \ldots, \dif{max} \, .
\end{equation*}
yields the sequence 
\begin{equation}\label{eq:Q_sequ}
    \begin{aligned}
        \Qc{A(0)} & \crksub{1} \Qc{A(1)} \crksub{1} \cdots 
        \crksub{1} \Qc{A(\delta)} \crksub{2} \Qc{A(\delta+1)} \crksub{2} 
        \cdots \\
        \cdots & \crksub{2} \Qc{A(p^3-k^3)} \crksub{1} \Qc{A(p^3-k^3+1)} 
        \crksub{1} \cdots \\
        \cdots & \crksub{1} \Qc{A(\dif{max})} = \Qc{R-1} = \Spanbracket{\D x},
    \end{aligned}
\end{equation}
relevant for the geometric characterization.
\begin{remark}
    Note that $\Qc{A(0)}=\Spanbracket{\D \varphi_{[0,K-1]}}$ and \mbox{$\Qc{A(\dif{max})}=\Qc{R-1}$} 
    are always integrable, because $\varphi_{[0,K-1]}(x)$ depends solely on the state $x$,
    whereas \mbox{$\Qc{R-1} = \Spanbracket{\D x}$} is directly implied by~\eqref{eq:imp_exterior_deriv}.
\end{remark}

Using~\eqref{eq:Q_sequ}, the following theorem provides a
geometric characterization of the structurally flat triangular 
form \eqref{eq:GTF_3_inputs}.
 
\begin{theorem}\label{thm:SFE_to_GTF3}
    A system~\eqref{eq:f_xu_affine_three_inputs} with a flat 
    output~\eqref{eq:x_flat_output_three_inputs} is SFE 
    to~\eqref{eq:GTF_3_inputs} if and only if the 
    sequence~\eqref{eq:Q_sequ} consists 
    purely of integrable codistributions.
\end{theorem}

Appendix~\ref{app:prf_thm_gtf3} provides a proof sketch of 
Theorem~\ref{thm:SFE_to_GTF3}. 

\subsection{ Equivalence to the Triangular Form by Prolongations }

Given a system for which 
Theorem~\ref{thm:SFE_to_GTF3} does not hold, there exist 
sufficient conditions to render a system SFE to~\eqref{eq:GTF_3_inputs}
by prolongations of the original control inputs.

\newpage
\begin{corollary}\label{cor:prolong}
    Let $y=(\varphi^1(x),\varphi^2(x),\varphi^3(x))$ be an 
    \mbox{$x$-flat} output of a control-affine system of the 
    form~\eqref{eq:f_xu_affine_three_inputs} that is not SFE 
    to~\eqref{eq:GTF_3_inputs}. If the
    given system is $\varphi$-SFL via minimal prolongations,
    then the system becomes SFE to~\eqref{eq:GTF_3_inputs}
    after at most $\ddiff{\varphi}$-fold prolongations of each input.
    Further, the components of the flat output $\varphi(x)$ 
    appear as $(z_1^1, z_2^1, z_3^1)$ in the corresponding
    triangular form defined by~\eqref{eq:GTF_3_inputs}.
\end{corollary}
\smallskip
\begin{remark}
    Corollary~\ref{cor:prolong} also applies to nonlinear 
    three-input systems that do not admit a control-affine 
    representation, except that at most $\ddiff{\varphi}+1$ 
    prolongations of each input are required, since a one-fold 
    prolongation of each input yields a control-affine system.
\end{remark}
Appendix~\ref{app:prf_thm_prolong} provides a sketch of the proof
for Corollary~\ref{cor:prolong}. 

% Whether $x$-flatness implies $\varphi$-SFL via
% minimal prolongations remains an open question.
% % \footnote{To the best 
% % of the authors' knowledge, no counterexample has been found to date.}
% Given that---to the best of the authors' knowledge---no counterexample has
% been found to date, the triangular form~\eqref{eq:GTF_3_inputs} is a 
% candidate for a universal structurally flat triangular form for $x$-flat
% three-input systems. For systems that are $\varphi$-SFL via minimal 
% prolongations, Corollary~\ref{cor:prolong} guarantees that SFE 
% to~\eqref{eq:GTF_3_inputs} can always be achieved after 
% finitely many input prolongations. 

As discussed in \cite{HartlLinearizationFlatMultiInput2026}, whether 
$x$-flatness implies $\varphi$-SFL via
minimal prolongations remains an open question, although---or perhaps 
because---no counterexample has been found to date.
Since systems that are $\varphi$-SFL via minimal 
prolongations always become SFE to~\eqref{eq:GTF_3_inputs} after 
finitely many input prolongations according to Corollary~\ref{cor:prolong},
the proposed triangular form is a 
candidate for a universal structurally flat triangular form for $x$-flat
three-input systems. 
This motivates the search for 
computationally tractable necessary and sufficient 
conditions characterizing~\eqref{eq:GTF_3_inputs} directly by the system dynamics, without 
requiring prior knowledge of a valid $x$-flat output.
\begin{remark}
    Such conditions would yield an iterative flatness test: verify SFE 
    to~\eqref{eq:GTF_3_inputs}, and if the test fails, extend the 
    system by a one-fold prolongation of each input and repeat. 
    Given the upper bound $\ddiff{\varphi} \leq 2n-6$ that can be 
    derived via, 
    e.g., \cite[Theorem~4]{kolar_properties_2016}, 
    this procedure terminates after at most $2n-6$ iterations for 
    any system that is $x$-flat and $\varphi$-SFL via minimal 
    prolongations. If the test fails at every step, the system is 
    either not $x$-flat or not $\varphi$-SFL via minimal 
    prolongations.
\end{remark}
\section{Examples}
\label{sec:examples}
 
We illustrate our results on two examples. The first is a 
physical system that is directly SFE the proposed triangular form, 
while the second is an academic example that requires input 
prolongations before SFE to~\eqref{eq:GTF_3_inputs} can be achieved.

%%%%%%%%%%%%%%%%%%%% EXAMPLE 1 %%%%%%%%%%%%%%%%%%%%%%%%%%%%%%%%%%
\smallskip
 \begin{example}[Planar Aerial Manipulator]%
\label{ex:aerial_manipulator}
Consider the planar aerial manipulator studied, e.g., 
in~\cite{welde_role_2023, 
HartlExactLinearizationMinimally2024a}, whose equations of 
motion are of the form~\eqref{eq:f_xu_affine_three_inputs} 
with $n=8$ states and three inputs. As shown 
in~\cite{welde_role_2023}, a flat output is given by 
\mbox{$y = (\theta,\; z_e - r_e\sin(\phi+\theta),\; 
x_e - r_e\cos(\phi+\theta))$}. For the given flat output 
we have \mbox{$K=(2,2,2)$}, \mbox{$R=(4,4,4)$} (see \cite{HartlExactLinearizationMinimally2024a}), 
yielding \mbox{$\ddiff{\varphi}=4$}.
Since \mbox{$R - K = (2,2,2)$}, the 
arrangement~\eqref{eq:rearrangement_case_rank1} is satisfied 
with \mbox{$\dif{max} = \dif{min} = 2$} and \mbox{$\delta = 0$}. After 
applying the input transformation 
\mbox{$\hat u^1 = \varphi^1_{[2]}(x,u)$}, the flat-output derivatives 
take the form~\eqref{eq:deriv_form_case_rank1} with 
\mbox{$P = (2,2,2) = K$} and $s = 2$. The multi-index 
$A(j)$ from~\eqref{eq:multi_index_A} evaluates to
\begin{equation*}
    A(0) = (1,1,1), \quad A(1) = (2,2,2), \quad 
    A(2) = (3,3,3),
\end{equation*}
yielding $d_{rp} = r^2 - p^2 = 2$ and confirming 
\mbox{$n = k^1 + p^2 + r^3 = 2 + 2 + 4 = 8$}. 

The codistribution sequence~\eqref{eq:Q_sequ} reduces to
\begin{equation*}
    \Qc{A(0)} \crksub{1} \Qc{A(1)} \crksub{1} 
    \Qc{A(2)} = \Spanbracket{\D x},
\end{equation*}
being completely integrable.  
By Theorem~\ref{thm:SFE_to_GTF3} the system is SFE 
to~\eqref{eq:GTF_3_inputs} yielding the triangular form
\begin{equation*}
    \begin{alignedat}{3}
        \dot z^1_1 & = z^2_1, & \quad 
        \dot z^1_2 & = z^2_2, & \quad 
        \dot z^1_3 & = z^2_3, \\
        \dot z^2_1 & = \hat u^1, & 
        \dot z^2_2 & = \hat u^2, & 
        \dot z^2_3 & = a^2_3(z_3^1, z_3^2, z_3^3)\!+\! b^2_{3,1}(z_3^1, z_3^3)\hat u^1 \\ 
        &&&&& \hspace{1em}  
        \!+\! b^2_{3,2}(z_3^1, z_3^3)\hat u^2, \\
        & & & & \dot z^3_3 & = z^4_3, \\
        & & & & \dot z^4_3 & = \hat u^3 ,
    \end{alignedat}
\end{equation*}
with subsystem dimensions $k^1 = 2$, $p^2 = 2$, and $r^3 = 4$.
% Note that since 
% $\delta = 0$, the $z_2$-subsystem consists only of the 
% integrator chain $\dot z^1_2 = z^2_2$, $\dot z^2_2 = \hat u^2$
% --- the coupled dynamics appear exclusively in the 
% $z_3$-subsystem.
\end{example}
%%%%%%%%%%%%%%%%%%%% END EXAMPLE 1 %%%%%%%%%%%%%%%%%%%%%%%%%%%%%%%%%%

\smallskip

%%%%%%%%%%%%%%%%%%%% EXAMPLE 2 %%%%%%%%%%%%%%%%%%%%%%%%%%%%%%%%%%
\begin{example}[Academic Example]\label{ex:academic}
Consider the system 
from~\cite[Sec.~4.4.2]{kolar_contributions_2017}:
\begin{equation}\label{eq:exa_n_7_m_3_system}
    \begin{alignedat}{2}
        \dot x^1 &= u^1, & \hspace{1.5em} 
        \dot x^5 &= -x^6 + x^4x^7u^1, \\
        \dot x^2 &= x^3+x^4u^1, & 
        \dot x^6 &= -x^5u^1 + x^7(u^1u^3 \\
        \dot x^3 &= u^2 - u^1u^3, & 
        & \;\;{-}\,u^2{-}1) + (x^4{+}u^1)x^4u^1, \\
        \dot x^4 &= u^3, & 
        \dot x^7 &= x^4+u^1.
    \end{alignedat}
\end{equation}
System~\eqref{eq:exa_n_7_m_3_system} admits the flat output 
$\varphi = (x^2, x^1, x^5)$ with $K=(1,1,1)$, $R=(4,3,4)$, 
and $\ddiff{\varphi} = 4$. Given \mbox{$R-K=(3,2,3)$}, 
the components are already arranged such 
that~\eqref{eq:rearrangement_case_rank1} holds. 
Since~\eqref{eq:exa_n_7_m_3_system} does not admit a 
control-affine representation\footnote{This means, there exists no invertible 
input transformation~\mbox{$\tilde{u} = \phi_{\tilde u}(x,u)$} such 
that~\eqref{eq:exa_n_7_m_3_system} takes a control-affine form.
See, e.g., \cite[Lemma~2.2]{gstottner_analysis_2023}, for a test to prove
whether a system $\dot x = f(x,u)$ allows a control-affine representation.}, 
we first perform a one-fold
prolongation of each input, yielding
\begin{equation}\label{eq:exa_prolonged}
    \dot x = f(x,u), \quad \dot u = u_{[1]},
\end{equation}
with extended state $x_e=(x,u)$ and input $u_{[1]}$. Note that
prolongation preserves the flat output 
$\varphi = (x^2, x^1, x^5)$ with $\ddiff{\varphi} = 4$, but 
now $K_e=(2,2,2)$ and $R_e=(5,4,5)$, giving $R_e-K_e=(3,2,3)$, 
$\dif{max} = 3$, $\dif{min} = 2$, and $\delta = 1$. The 
multi-index $A(j)$ from~\eqref{eq:multi_index_A} evaluates to
\mbox{$A(0) = (1,1,1)$}, \mbox{$A(1) = (2,1,2)$},
\mbox{$A(2) = (3,2,3)$}, and \mbox{$A(3) = (4,3,4)$}.
Then the codistribution sequence~\eqref{eq:Q_sequ} takes the form
\begin{equation*}
    \Qc{A(0)} \crksub{1} \Qc{A(1)} \crksub{2} 
    \Qc{A(2)} \crksub{1} \Qc{A(3)} 
    = \Spanbracket{\D x, \D u},
\end{equation*}
where
\begin{equation}\label{eq:Q_seq_academic}
    \begin{aligned}
        \Qc{A(0)} &= \langle \D x^1, \D x^2, \D x^5, 
        \D x^3 {+} u^1\D x^4, \\[-0.5ex]
        & \hspace{3em} {-}x^7\D x^3 {-} \D x^6 
        {+} x^4u^1\D x^7, \D u^1 \rangle, \\
        \Qc{A(1)} &= \Qc{A(0)} + \langle \D x^7 \rangle, \\
        \Qc{A(2)} &= \Qc{A(1)} + \langle \D u^2, 
        \D x^6 \rangle, \\
        \Qc{A(3)} &= \Qc{A(2)} + \langle \D u^3 \rangle.
    \end{aligned}
\end{equation}
Given the state $x_e=(x,u)$ of the extended system \eqref{eq:exa_prolonged},
it follows that all codistributions are integrable. According to
Theorem~\ref{thm:SFE_to_GTF3}, the prolonged system is SFE 
to~\eqref{eq:GTF_3_inputs}. To construct the transformation, 
note that $(z_j^1, z_j^2) = \varphi^j_{[0,1]}$ for 
$j \in \{1,2,3\}$ and the remaining coordinates are determined 
by the growing sequence~\eqref{eq:Q_seq_academic}:
\begin{equation}\label{eq:trf_academic}
    \begin{alignedat}{3}
        z_1^1 & = x^2, & \hspace{1em}
        z_2^1 & = x^1, & \hspace{1em}
        z_3^1 & = x^5, \\
        z_1^2 & = x^3{+}x^4u^1, & 
        z_2^2 & = u^1, & 
        z_3^2 & = x^4x^7u^1{-}x^6, \\[0.5ex]
        & & & & 
        z_3^3 & = x^7, \\
        & & 
        z_2^3 & = u^2, & 
        z_3^4 & = x^6, \\
        & & & & 
        z_3^5 & = u^3.
    \end{alignedat}
\end{equation}
Applying~\eqref{eq:trf_academic} together with the input 
transformation 
\begin{equation*}
    (\hat u^1, \hat u^2, \hat u^3) 
    = (\varphi^1_{[2]}(x_e, u_{[1]}), u^2_{[1]}, u^3_{[1]})\, ,
\end{equation*}
the $z_1$-subsystem takes the form 
$\dot z_1^1 = z_1^2$, $\dot z_1^2 = \hat u^1$, and the 
$z_2$- and $z_3$-subsystems read
\begin{equation}\label{eq:gtf3_academic}
    \begin{alignedat}{2}
        \dot z_2^1 & = z_2^2, & \hspace{1em}
        \dot z_3^1 & = z_3^2, \\
        & & 
        \dot z_3^2 & = z_2^2z_3^1 {+} z_3^3 
        {+} z_3^3\hat{u}^1, \\
        \dot z_2^2 & = a_2^2(\overline{z}_2^{2,3}, 
        \overline{z}_3^{2,4}) & 
        \dot z_3^3 & = a_3^3(z_2^2, 
        \overline{z}_3^{2,4}), \\
        & \hspace{1em} {+}\, b_2^2(z_2^2, 
        \overline{z}_3^{2,4})\hat{u}^1, & & \\
        \dot z_2^3 & = \hat{u}^2, &
        \dot z_3^4 & = a_3^4(\overline{z}_2^{2,3}, 
        \overline{z}_3^{2,5}), \\
        & &
        \dot z_3^5 & = \hat{u}^3,
    \end{alignedat}
\end{equation}
where $\overline{z}_2^{2,3} = (z_2^{2}, z_2^{3})$ and
$\overline{z}_3^{2,4} = (z_3^{2}, z_3^{3}, z_3^{4})$. The 
system~\eqref{eq:gtf3_academic} is of the 
form~\eqref{eq:GTF_3_inputs} with $k^1=2$, $p^2=3$, and 
$r^3=5$. Note that $\delta = 1$ is reflected by the way
the flat parameterization \eqref{eq:flat_parametrization} is derived
from top to bottom: 
starting from the integrator chains, one first obtains the flat parameterization
for $(z_1^1, z_1^2)$, $(z_2^1, z_2^2)$, $(z_3^1, z_3^2)$, and 
$\hat u^1$. At the next level, only 
$z_3^3 = z_3^{k^3+\delta}$ can be recovered from the dynamics 
of $\dot z_3^2$, before $z_2^3 = z_2^{k^2+1}$ and 
$z_3^4 = z_3^{k^3+\delta+1}$ follow from $\dot z_2^2$ and $\dot z_3^3$.
\end{example}
%%%%%%%%%%%%%%%%%%%% END EXAMPLE 2 %%%%%%%%%%%%%%%%%%%%%%%%%%%%%%%%%%

%% ================================================================
\section{Conclusion}
\label{sec:conclusion}
 
We introduced a general structurally flat triangular form for 
\mbox{$x$-flat} control-affine three-input systems and derived 
necessary and sufficient conditions for static 
feedback equivalence to the proposed form that are
formulated via the 
integrability of a sequence of codistributions based on a 
given flat output. We showed that every 
\mbox{$x$-flat} three-input system that is $\varphi$-SFL via 
minimal prolongations can be brought into the general triangular 
form after a finite number of input prolongations. The results 
were illustrated on a mechanical and an academic example.
For the two-input case, the analogous triangular 
form~\cite{gstottner_triangular_2024} has been used to develop
distribution-based algorithms to identify flat-output 
candidates~\cite{HartlTriangularFormsXFlat2025}.
Future work will therefore focus on extending these results to the
three-input case using~\eqref{eq:GTF_3_inputs}.

\bibliography{IEEEabrv,mybibfile}

\newpage
%% ================================================================
\appendices

\section{Proof of Theorem~\ref{thm:SFE_to_GTF3}}
\label{app:prf_thm_gtf3}

\emph{Necessity:} It suffices to verify that the 
sequence~\eqref{eq:Q_sequ} for any system of 
the form~\eqref{eq:GTF_3_inputs} with flat output 
\mbox{$\varphi = (z_1^1, z_2^1, z_3^1)$} consists 
solely of integrable codistributions.

\emph{Sufficiency:} We show that any 
system~\eqref{eq:f_xu_affine_three_inputs} with a flat 
output~\eqref{eq:x_flat_output_three_inputs} for which all 
codistributions in~\eqref{eq:Q_sequ} are integrable can be 
brought into the form~\eqref{eq:GTF_3_inputs}.

\textbf{Step~1} (State transformation):
Since the codistributions in~\eqref{eq:Q_sequ} are integrable 
by assumption, there exist functions $g_2^i(x)$, 
$i = k^2{+}1, \ldots, p^2$, and $g_3^i(x)$, 
$i = k^3{+}1, \ldots, r^3$, such that
\begin{equation*}
\begin{aligned}
    \Qc{A(1)} &= \Qc{A(0)} + \langle \D g^{k^3+1}_3 
      \rangle, \\[-0.6em]
    &\hspace{0.5em} \vdots \\[-0.6em]
    \Qc{A(\delta)} &= \Qc{A(\delta-1)} + \langle 
      \D g^{k^3+\delta}_3 \rangle, \\
    \Qc{A(\delta+1)} &= \Qc{A(\delta)} + \langle 
      \D g^{k^2+1}_2, \D g^{k^3+\delta+1}_3 \rangle, \\[-0.6em]
    &\hspace{0.5em} \vdots \\[-0.6em]
    \Qc{A(p^3\!-\!k^3)} &= \Qc{A(p^3\!-\!k^3\!-\!1)} 
      + \langle \D g^{p^2}_2, \D g^{p^3}_3 \rangle, \\
    \Qc{A(p^3\!-\!k^3\!+\!1)} &= \Qc{A(p^3\!-\!k^3)} 
      + \langle \D g^{p^3+1}_3 \rangle, \\[-0.6em]
    &\hspace{0.5em} \vdots \\[-0.6em]
    \Qc{A(\dif{max})} &= \Qc{A(r^3\!-\!k^3\!-\!1)} 
      + \langle \D g^{r^3}_3 \rangle 
      = \langle \D x \rangle.
\end{aligned}
\end{equation*}
We introduce the invertible state transformation
\begin{equation}\label{eq:state_trf_x_to_z}
\begin{aligned}
    (z_1^1, \ldots, z_1^{k^1}, z_2^1, \ldots, z_2^{k^2}, 
    z_3^1, \ldots, z_3^{k^3}) &= \varphi_{[0,K-1]}(x), \\
    z_2^{k^2+1} = g_2^{k^2+1}(x), \ldots, z_2^{p^2} 
    &= g_2^{p^2}(x), \\
    z_3^{k^3+1} = g_3^{k^3+1}(x), \ldots, z_3^{r^3} 
    &= g_3^{r^3}(x).
\end{aligned}
\end{equation}
In these coordinates, the $\Qc{}$-sequence takes the form
\begin{equation*}\label{eq:Q_sequence_proof}
\begin{aligned}
    % \Qc{A(0)} &= \langle \D \overline{z}_1^{k^1}, 
    %   \D \overline{z}_2^{k^2}, 
    %   \D \overline{z}_3^{k^3} \rangle, \\
    \Qc{A(1)} &= \langle \D z_1, 
      \D \overline{z}_2^{k^2}, 
      \D \overline{z}_3^{k^3+1} \rangle, \\[-0.6em]
    &\hspace{0.5em} \vdots \\[-0.6em]
    % \Qc{A(\delta)} &= \langle \D \overline{z}_1^{k^1}, 
    %   \D \overline{z}_2^{k^2}, 
    %   \D \overline{z}_3^{k^3+\delta} \rangle, \\
    \Qc{A(\delta+1)} &= \langle \D z_1, 
      \D \overline{z}_2^{k^2+1}, 
      \D \overline{z}_3^{k^3+\delta+1} \rangle, \\[-0.6em]
    &\hspace{0.5em} \vdots \\[-0.6em]
    \Qc{A(p^3\!-\!k^3)} &= \langle \D z_1, 
      \D z_2, 
      \D \overline{z}_3^{p^3} \rangle, \\[-0.6em]
    &\hspace{0.5em} \vdots \\[-0.6em]
    \Qc{A(\dif{max})} &= \langle \D z_1, 
      \D z_2, 
      \D z_3 \rangle.
\end{aligned}
\end{equation*}

%%%%%%%%%%%%%%%%%%%%%%%%%%%%%%%%%%%%%%%%%%%%%%%%%%%%%%%%%%%%%%%%%%%%%%%%%%%%%%%%%
\begin{figure*}[!b]
\hrulefill
\normalsize
\begin{equation}
    \tag{$DS_{3e}$}
    \label{eq:deriv_form_extended_s_drp}
    \begin{alignedat}{3}
        % R+1
        \varphi^1_{[r^1\!+\!1]} & = \hat u^1_{[d_{\max}\!+\!1]} & 
        \hspace{1em} & \varphi^2_{[r^2\!+\!1]}(x, \hat u^1_{[0,d_{\min}\!+\!1]},\, u^2_{[0,d_{rp}\!+\!1]},\, u^3_{[0,1]}) & 
        \hspace{1em} & \varphi^{3}_{[r^3\!+\!1]}( x, \hat u^1_{[0,d_{\max}\!+\!1]},\, u^2_{[0,d_{rp}\!+\!1]},\, u^3_{[0,1]}) \\[-1ex]
        % dots
        &  && \hspace{0.5em} \vdots &&  \\[-1ex]
        % intermediate step
        &\hspace{0.5em} \vdots&
        & \varphi^2_{[r^2\!+\!d_{rp}]}(x, 
        \hat u^1_{[0,r^2\!-\!k^2\!+\!d_{rp}]},\, 
        u^2_{[0,2d_{rp}]},\, u^3_{[0,d_{rp}]}) &
         & \hspace{0.5em} \vdots \\[-1ex]
        % dots
        &  && \hspace{0.5em} \vdots &&  \\[-1ex]
        % K_e
        \varphi^1_{[k^1\!+\!d]} & = \hat u^1_{[d]} & 
        & \varphi^2_{[k^2\!+\!d]}(x, \hat u^1_{[0,d]},\, u^2_{[0,2d_{rp}\!+\!\delta]},\, u^3_{[0,d_{rp}\!+\!\delta]}) & 
        & \varphi^3_{[k^3\!+\!d]}( x, \hat u^1_{[0,d]},\, u^2_{[0,2d_{rp}]},\, u^3_{[0,d_{rp}]}) \\[-1ex]
        % dots
        & \hspace{0.5em} \vdots && && \hspace{0.5em} \vdots \\[-1ex]
        % \delta
        \varphi^1_{[k^1\!+\!d+\!\delta]} & = \hat u^1_{[d\!+\!\delta]} & 
        & \hspace{0.5em} \vdots & 
        & \varphi^3_{[k^3\!+\!d\!+\!\delta]}( x, \hat u^1_{[0,d\!+\!\delta]},\, u^2_{[0,2d_{rp}\!+\!\delta]},\, u^3_{[0,d_{rp}\!+\!\delta]}) \\[-1ex]
        % dots
        & \hspace{0.5em} \vdots &&  && \hspace{0.5em} \vdots \\[-1ex]
        % R_e-1
        \varphi^1_{[r^1\!+\!d\!-\!1]} & = \hat u^1_{[d_{\max}\!+\!d\!-\!1]} & 
        & \varphi^2_{[r^2\!+\!d\!-\!1]}(x, \hat u^1_{[0,d_{\min}\!+\!d\!-\!1]},\, u^2_{[0,d_{rp}\!+\!d\!-\!1]},\, u^3_{[0,d\!-\!1]}) & 
        & \varphi^3_{[r^3\!+\!d\!-\!1]}( x, \hat u^1_{[0,d_{\max}\!+\!d\!-\!1]},\, u^2_{[0,d_{rp}\!+\!d\!-\!1]},\, u^3_{[0,d\!-\!1]})
    \end{alignedat}
\end{equation}
\end{figure*}
%%%%%%%%%%%%%%%%%%%%%%%%%%%%%%%%%%%%%%%%%%%%%%%%%%%%%%%%%%%%%%%%%%%%%%%%%%%%%%%%%

\textbf{Step~2} (Dynamics for 
$z_2^{k^2},\!\ldots,\!z_2^{p^2\!-\!1}$ and
$z_3^{k^3},\!\ldots,\!z_3^{p^3\!-\!1}$ ):
We apply the static input transformation 
$\hat u^1 = \varphi^1_{[k^1]}(z,u)$ replacing~$u^1$ after a possible
relabeling.
For convenience, we denote the vector 
field~\eqref{eq:extended_vector_field} by $f_u$ regardless of 
the choice of state and input coordinates. In the $z$-coordinates, 
$f_u$ takes the form
\begin{equation}\label{eq:vec_field_Step_2}
\begin{aligned}
    f_{u} &= z_j^2\pad{z_j^1} + \cdots 
      + z_j^{k^j}\pad{z_j^{k^j-1}} 
      + f_l^{i_l}(z, \hat u^1, u^2, u^3)\pad{z_l^{i_l}} \\
    &\quad + \sum_{\alpha=0}^{\dif{max}-1} 
      \hat u^1_{[\alpha+1]}\pad{\hat u^1_{[\alpha]}} 
      + \sum_{\alpha=0}^{d_{rp}-1} 
      u^l_{[\alpha+1]}\pad{u^l_{[\alpha]}} \, ,
\end{aligned}
\end{equation}
with $j=1,2,3$, $l=2,3$, $i_2 = k^2, \ldots, p^2$, and 
\mbox{$i_3 = k^3, \ldots, r^3$}. 

The key idea underlying the remainder of this proof is 
as follows: since $\Qc{A(i)} \subset \Pc{A(i)}$, 
each $\Pc{A(i)}$ admits an explicit basis. 
By comparing this basis with the relation 
$\D \varphi^j_{[i]} = \D \Lie_{f_u}\varphi^j_{[i-1]}$, 
one can determine the admissible dependencies of 
the coefficient functions of~\eqref{eq:vec_field_Step_2}.
For each $i \in \{0, \ldots, \delta-1\}$, 
a basis for the codistribution
\begin{equation}\label{eq:Pc_Aj_delta}
    \Pc{A(i+1)} = \Pc{A(i)} 
    + \Spanbracket{\D \varphi^1_{[k^1+i]}, 
      \D \varphi^3_{[k^3+i]}}
\end{equation}
is given by
\begin{equation}\label{eq:Pc_Aj_delta_z}
    \Pc{A(i+1)} = \langle \D z_1, 
      \D \hat u^1_{[0,i]}, 
      \D \overline{z}_2^{k^2}, 
      \D \overline{z}_3^{k^3+i+1} \rangle.
\end{equation}

Given that $\varphi^3_{[k^3+i]} 
= \Lie_{f_u}\varphi^3_{[k^3+i-1]}$, it follows that
$f_3^{k^3+i}$ must explicitly depend on 
$z_3^{k^3+i+1}$ while being independent of 
$(z_2^{k^2+1}, \ldots, z_2^{p^2})$ and 
$(z_3^{k^3+i+2}, \ldots, z_3^{r^3})$. This proves 
the regularity condition \eqref{eq:reg_cond_1}.

For $i \in \{0, \ldots, p^2-k^2-1\}$, the codistributions
\begin{equation*}\label{eq:Pc_Ai_coupled}
    \Pc{A(\delta+i+1)} = \Pc{A(\delta+i)} 
      + \Spanbracket{
      \D \varphi^1_{[k^1+\delta+i]}, 
      \D \varphi^2_{[k^2+i]}, 
      \D \varphi^3_{[k^3+\delta+i]} }
\end{equation*}
take the form
\begin{equation*}\label{eq:P_seq_coupled}
    \Pc{A(\delta+i+1)} = \langle \D z_1, 
      \D \hat u^1_{[0,\delta+i]}, 
      \D \overline{z}_2^{k^2+i+1}, 
      \D \overline{z}_3^{k^3+\delta+i+1} \rangle.
\end{equation*}

Again, given that $\varphi^2_{[k^2+i]} 
= \Lie_{f_u}\varphi^2_{[k^2+i-1]}$ and
$\varphi^3_{[k^3+\delta+i]} 
= \Lie_{f_u}\varphi^3_{[k^3+\delta+i-1]}$, it follows that
$f_2^{k^2+i}$ and $f_3^{k^3+\delta+i}$ must explicitly 
depend on $(z_2^{k^2+i+1},z_3^{k^3+\delta+i+1})$ such that
the regularity condition \eqref{eq:reg_cond_2} is satisfied.

It remains to determine the dependency of 
$f_3^{k^3+i}$, \mbox{$i\in\{0,\ldots,p^3-k^3-1\}$}, on the 
$z_2$-subsystem. Given the $z_2$-integrator chain
\mbox{$\dot z_2^l = z_2^{l+1}$} for \mbox{$1 \leq l \leq k^2-1$}, 
any dependency of 
$f_3^{k^3+i}$ on $z_2^l$ propagates through successive Lie 
derivatives and after $k^2 - l$ differentiations, 
it reaches $z_2^{k^2}$. One further differentiation introduces 
$f_2^{k^2}$, which depends on $z_2^{k^2+1}$ and/or 
$z_3^{k^3+\delta+1}$. Since 
$\D z_2^{k^2+1}$ and $\D z_3^{k^3+\delta+1}$
first appear in 
$\Pc{A(\delta+1)}$, they must not enter
$\Pc{A(i)}$ for $i \leq \delta$, implying that
$f_3^{k^3+\delta-1}$ may depend on at most $\overline{z}_2^{k^2}$.
This requires $l+\delta-1-i \leq k^2$, which further yields that 
$f_3^{k^3+i}$ may depend on at most $\overline{z}_2^{k^2-\delta+1+i}$, for
\mbox{$i\in\{0,\ldots,p^3-k^3-1\}$}.

Given the control-affine structure of the system, 
the functions
$f_3^{k^3+i}$, $i \in \{0, \ldots, p^3-k^3-1\}$, 
decompose as
\begin{equation*}
\begin{aligned}
    f_3^{k^3+i} &= a_3^{k^3+i}(z_1, 
      \overline{z}_2^{k^2-\delta+i+1}, 
      \overline{z}_3^{k^3+i+1}) 
      + b_{3,1}^{k^3+i}(\ldots)\hat u^1,
\end{aligned}
\end{equation*}
while $f_2^{k^2+i}$, $i \in \{0, \ldots, p^2-k^2-1\}$ 
take the form
\begin{equation*}
\begin{aligned}
    f_2^{k^2+i} &= a_2^{k^2+i}(z_1, 
      \overline{z}_2^{k^2+i+1}, 
      \overline{z}_3^{k^3+\delta+i+1}) 
      + b_{2,1}^{k^2+i}( \ldots )\hat u^1 .
\end{aligned}
\end{equation*}

\textbf{Step~3} (Input transformation 
$\hat u^2 = f^{p^2}_2(z,\hat u^1, u^2, u^3)$): 

Given~\eqref{eq:deriv_form_case_rank1}, 
$\varphi^j_{[p^j]}$, \mbox{$j\in \{2,3\}$}, are the first functions 
depending explicitly on $u^2$ and additionally on $u^3$ if $s=0$. 
Since the codistribution
\mbox{$\Pc{A(p^3-k^3+1)} = 
\Spanbracket{\D \varphi^1_{[0,p^1]}, 
\D \varphi^2_{[0,p^2]}, \D \varphi^3_{[0,p^3]}}$}
is spanned by
\[
\Pc{A(p^3-k^3+1)} \!=\! \Spanbracket{\D z_1,
\D \hat u^1_{[0,p^3-k^3]}, 
\D z_2, \alpha \D u^2 
+ \beta \D u^3, \D \overline{z}_3^{p^3+1}},\]
with $\beta\!=\!0$ for $s\!>\!0$, the functions $f_2^{p^2}$ and 
$f_3^{p^3}$ are given by
\begin{equation*}\label{eq:fct_f_step_p2mk2}
\begin{aligned}
    f_j^{p^j} &= \tilde{a}_j(z_1, z_2, 
      \overline{z}_3^{p^3+1}) 
      + \tilde{b}_{j,1}( \ldots )\hat u^1 \\
      & \hspace{2em}
      + \tilde{b}_{j,2}(\ldots) u^2 
      + \tilde{b}_{j,3}( \ldots ) u^3,
\end{aligned}
\end{equation*}
$j \in \{2,3\}$. Applying the static input transformation
\begin{equation}\label{eq:input_trf_Step3}
    \hat u^2 = f_2^{p^2}(z_1, z_2, 
      \overline{z}_3^{p^3+1}, \hat u^1, u^2, u^3),
\end{equation}
replacing $u^2$ yields 
$f_2^{p^2} = \hat u^2$ and
\begin{equation*}
    f_3^{p^3} = a_3^{p^3}(z_1, z_2, 
      \overline{z}_3^{p^3+1}) 
      + b_{3,1}^{p^3}(\ldots)\hat u^1 
      + b_{3,2}^{p^3}(\ldots)\hat u^2.
\end{equation*}
Note that $f_3^{p^3}$ is independent of $u^3$, as otherwise 
$\D z_3^{p^3+1}$ would not be contained 
in~$\Pc{A(p^3\!-\!k^3\!+\!1)}$.

\textbf{Step~4} (Remaining dynamics for 
$z_3^{p^3+1}, \ldots, z_3^{r^3}$):
After the input transformation~\eqref{eq:input_trf_Step3}, 
the subsequent codistributions grow by one $z_3$-coordinate 
at each level:
\begin{equation*}
\begin{aligned}
    \Pc{A(p^3\!-\!k^3\!+\!1)} &= \Pc{A(p^3\!-\!k^3)} 
      + \langle \D \hat u^1_{[p^3-k^3]}, 
      \D \hat u^2,
      \D z_3^{p^3+1} \rangle, \\[-0.6em]
    % \Pc{A(p^3\!-\!k^3\!+\!2)} &= 
    %   \Pc{A(p^3\!-\!k^3\!+\!1)} + \langle 
    %   \D \hat u^1_{[p^3-k^3+1]}, 
    %   \D \hat u^2_{[1]}, 
    %   \D z_3^{p^3+2} \rangle, \\[-0.6em]
    &\hspace{0.6em} \vdots \\[-0.6em]
    \Pc{A(\dif{max}+1)} &= 
      \Pc{A(\dif{max})} + \langle 
      \D \hat u^1_{[\dif{max}]}, 
      \D \hat u^2_{[d_{rp}]}, \D u^3 \rangle.
\end{aligned}
\end{equation*}
Note that $\D u^3$ must appear in $\Pc{A(\dif{max}+1)}$, 
since flatness implies $\langle\D x, \D u \rangle \subset 
\langle\D \varphi_{[0,R]} \rangle$. Analogous reasoning as in 
Steps~2 and~3 shows that for 
$i \in \{p^3-k^3, \ldots, r^3-k^3-1\}$, the functions 
$f_3^{k^3+i}$ take the form
\begin{equation*}
\begin{aligned}
    f_3^{k^3+i} &= a_3^{k^3+i}(z_1, z_2, 
      \overline{z}_3^{k^3+i+1}) \\
      & \hspace{2em} + b_{3,1}^{k^3+i}( \ldots )\hat u^1 
      + b_{3,2}^{k^3+i}( \ldots )\hat u^2,
\end{aligned}
\end{equation*}
satisfying regularity condition~\eqref{eq:reg_cond_3}. Normalizing
$f_3^{r^3} = \hat u^3$ by a static input transformation
yields exactly the dynamics given in~\eqref{eq:GTF_3_inputs}, 
completing the proof.

\section{Proof of Corollary~\ref{cor:prolong}}
\label{app:prf_thm_prolong}

By assumption, we consider a control-affine three-input system of the form
\begin{equation*}
    \dot x = f(x) + g_1(x)\hat{u}^1 + g_2(x)u^2 + g_3(x)u^3
\end{equation*}
with an $x$-flat output
\mbox{$\varphi(x)=(\varphi^1(x),\varphi^2(x),\varphi^3(x))$}
that is $\varphi$-SFL via minimal prolongations. 
To the given flat output, we associate the multi-indices $K$ and $R$, and the 
differential difference $\ddiff{\varphi}$, which we briefly denote by $d$. 
According to Theorem~\ref{thm:rank1_case}, the components of $\varphi$ 
can be relabeled such that~\eqref{eq:rearrangement_case_rank1} holds
and the derivatives $\varphi_{[0,R]}$ can be brought into the 
form~\eqref{eq:deriv_form_case_rank1} with the integers $p^2$ and $p^3$. 
For simplicity, we assume the derivatives $\varphi_{[0,R]}$ are
already of the form~\eqref{eq:deriv_form_case_rank1} with $s=d_{rp}$,
i.e., that means that item \textit{c)} of 
Theorem~\ref{thm:rank1_case} is already satisfied with $u^2$.

With the multi-indices $A(j)$, $j=0,\ldots,r^3-k^3+1$, 
see \eqref{eq:multi_index_A}, one can construct the associated $\Pc{}$- 
and $\Qc{}$-sequence \eqref{eq:P_sequ} and \eqref{eq:Q_sequ}.
In the following, we show that a $d$-fold prolongation of each input
of the considered system ensures integrability of the 
$\Qc{}$-sequence~\eqref{eq:Q_sequ} of the extended system
\begin{equation}\label{eq:extended_sys_prolonged_proof}
    \begin{aligned}
        \dot x & = f(x) + g_1(x)\hat{u}^1 + g_2(x)u^2 + g_3(x)u^3, \\
        \dot{\hat{u}}^1_{[0,d-1]} & = \hat{u}^1_{[1,d]}, \hspace{0.5em}
        \dot u^2_{[0,d-1]} = u^2_{[1,d]}, \hspace{0.5em}
        \dot u^3_{[0,d-1]} = u^3_{[1,d]} \, ,
    \end{aligned}
\end{equation}
with state $x_e = (x, u_{[0,d-1]})$ and input $u_{[d]}$,
where \mbox{$u=(\hat u^1, u^2, u^3)$}.
Note that the prolonged system admits the same flat output $\varphi(x)$ 
with differential difference $d$. However, the corresponding multi-indices 
$K_e = K+d$, $R_e = R+d$, and $A_e(j) = A(j)+d$, for $j=0,\ldots,r^3-k^3+1$,
are increased by $d$. 

The structure of the time derivatives $\varphi_{[0,R_e]}$ follows 
from~\eqref{eq:deriv_form_case_rank1} with $s = d_{rp}$ and 
its extension~\eqref{eq:deriv_form_extended_s_drp}.
By means of $\varphi_{[0,R_e]}$ and the fact that flatness implies 
\mbox{$\Spanbracket{\D x, \D u} \subset \Spanbracket{\D \varphi_{[0,R]}}$},
it can be deduced that the $\Pc{}$-sequence associated with the extended 
system~\eqref{eq:extended_sys_prolonged_proof} takes the form
\begin{equation}\label{eq:P_seq_prolonged_proof}
    \begin{aligned}
        % \Pc{A_e(0)} &= \Spanbracket{\D \varphi_{[0,K_e\!-\!1]}(x_e)} \, , \\
        \Pc{A_e(1)} &= \Spanbracket{\D \varphi_{[0,K_e\!-\!1]}(x_e), 
        \D \hat u^1_{[d]}, \D u^3_{[d_{rp}]}} \, , \\
        & \hspace{0.5em} \vdots \\[-0.6ex]
        \Pc{A_e(\delta\!-\!1)} &= \Spanbracket{\D \varphi_{[0,K_e\!-\!1]}(x_e), 
        \D \hat u^1_{[d,d\!+\!\delta\!-\!2]},\, \D u^3_{[d_{rp},d_{rp}\!+\!\delta\!-\!2]}} \, , \\
        \Pc{A_e(\delta)} &= \Spanbracket{\D x, 
        \D \hat u^1_{[0,d\!+\!\delta\!-\!1]}, 
        \D u^2_{[0,2d_{rp}\!+\!\delta\!-\!1]}, 
        \D u^3_{[0,d_{rp}\!+\!\delta\!-\!1]}} \, , \\
        \Pc{A_e(\delta\!+\!1)} &= \Spanbracket{\D x, 
        \D \hat u^1_{[0,d\!+\!\delta]}, 
        \D u^2_{[0,2d_{rp}\!+\!\delta]}, 
        \D u^3_{[0,d_{rp}\!+\!\delta]}} \, , \\
        & \hspace{0.5em} \vdots \\[-0.6ex]
        \Pc{A_e(p^3\!-\!k^3\!+\!1)} &= \Spanbracket{\D x, 
        \D \hat u^1_{[0,2\dif{max}]}, 
        \D u^2_{[0,d]}, 
        \D u^3_{[0,\dif{max}]}} \, , \\
        & \hspace{0.5em} \vdots \\[-0.6ex]
        \Pc{A_e(r^3\!-\!k^3)} &= \Spanbracket{\D x, 
        \D \hat u^1_{[0,d\!+\!\dif{max}\!-\!1]}, 
        \D u^2_{[0,d\!+\!\dif{min}]}, 
        \D u^3_{[0,d\!-\!1]}} \, .
    \end{aligned}
\end{equation}

Intersecting the $\Pc{}$-sequence \eqref{eq:P_seq_prolonged_proof} with 
\mbox{$\langle \D x_{e} \rangle=\Spanbracket{\D x, \D \hat u^1_{[0,d-1]}, \D u^2_{[0,d-1]}, \D u^3_{[0,d-1]}}$} yields the associated $\Qc{}$-sequence 
\eqref{eq:Q_sequ}, where each codistribution is solely spanned by 
exact differentials.

\end{document}